\documentclass{amsart}%
\usepackage{amsfonts}
\usepackage{amsmath}
\usepackage{amssymb}
\usepackage{graphicx}%
\providecommand{\U}[1]{\protect\rule{.1in}{.1in}}
\newtheorem{theorem}{Theorem}
\theoremstyle{plain}

\newtheorem{definition}{Definition}
\newtheorem{example}{Example}

\newtheorem{lemma}{Lemma}

\newtheorem{remark}{Remark}

\numberwithin{equation}{section}
\begin{document}
\title[ ]{Nonlinear versions of Korovkin's Abstract Theorems}
\author{Sorin G. Gal}
\address{Department of Mathematics and Computer Science\\
University of Oradea\\
University\ Street No. 1, Oradea, 410087, Romania}
\email{galso@uoradea.ro, galsorin23@gmail.com}
\author{Constantin P. Niculescu}
\address{Department of Mathematics, University of Craiova\\
Craiova 200585, Romania}
\email{constantin.p.pniculescu@gmail.com}
\date{March 2, 2021}
\subjclass[2000]{41A35, 41A36, 41A63}
\keywords{Choquet' integral, Korovkin theorem, weak additivity, monotone operator,
sublinear operator, (locally-) compact metric space}

\begin{abstract}
In this paper we prove Korovkin type theorems for sequences of sublinear,
monotone and weak additive operators acting on function spaces $C(X),$ where
$X$ is a compact or a locally compact metric space. Our results are
illustrated by a series of examples.

\end{abstract}
\maketitle

\section{Introduction}

The celebrated theorem of Korovkin \cite{Ko1953}, \cite{Ko1960} gives us
conditions for uniform approximation of continuous functions on a compact
interval via sequences of positive linear operators. Precisely, if
$(T_{n})_{n}$ is a sequence of positive linear operators that map $C\left(
[0,1]\right)  $ into itself such that the sequence $(T_{n}(f))_{n}$ converges
to $f$ uniformly on $[0,1]$ for the three special functions $e_{k}%
:x\rightarrow x^{k}$, where $k=0,1,2,$ then this sequence also converges to
$f$ uniformly on $[0,1]$ for every $f\in C([0,1])$. This statement remains
true by replacing $C\left(  [0,1]\right)  $ by the space $C_{2\pi}%
(\mathbb{R})$ of continuous and 2$\pi$-periodic functions defined on
$\mathbb{R}$ and considering as test functions the triplet $1$, $\cos$ and
$\sin$.

Over the years, many generalizations of Korovkin theorem appeared, all in the
framework of linear functional analysis. A nice account on the present state
of art is offered by the authoritative monograph of F. Altomare and M. Campiti
\cite{AC1994} and the excellent survey of F. Altomare \cite{Alt2010}. See also
\cite{BucurPaltin}.

Inspired by the Choquet theory of integrability with respect to a nonadditive
measure, we proved in \cite{Gal-Nic-1} an extension of Korovkin's theory to a
class of sublinear operators, called by us Choquet type operators. The aim of
the present paper is to show that our results still work in a context
considerably more general.

We shall need some rudiments of ordered Banach theory that can be covered from
the textbook of Meyer-Nieberg \cite{MN}, or from the papers \cite{Gal-Nic-3}
and \cite{NO2020}.

Given a metric space $X,$ we will denote by $\mathcal{F}(X)$ the vector
lattice of all real-valued functions defined on $X,$ endowed with the
pointwise ordering. Some important vector sublattices of it are
\begin{align*}
C(X)  &  =\left\{  f\in\mathcal{F}(X):\text{ }f\text{ continuous}\right\}  ,\\
C_{b}(X)  &  =\left\{  f\in\mathcal{F}(X):\text{ }f\text{ continuous and
bounded}\right\}
\end{align*}
and $\mathcal{U}C_{b}(X)=\left\{  f\in\mathcal{F}(X):\text{ }f\text{ uniformly
continuous and bounded}\right\}  $. Notice that the spaces $C(X),C_{b}(X)~$and
$\mathcal{U}C_{b}(X)$ coincide when $X$ is a compact metric space. If $d$
denotes the metric on $X,$ then an important family of Lipschitz continuous
functions in $C(X)$ is%
\[
d_{x}:X\rightarrow\mathbb{R},\text{\quad}d_{x}(y)=d(x,y)~(x\in X).
\]
It is also worth noticing that the spaces $C_{b}(X)$ and $\mathcal{U}C_{b}(X)$
are Banach lattices with respect to the pointwise ordering and the sup norm,%
\[
\left\Vert f\right\Vert _{\infty}=\sup\left\{  \left\vert f(x)\right\vert
:x\in X\right\}  .
\]

See \cite{CN2014}. As is well known, all norms on the $N$-dimensional real
vector space $\mathbb{R}^{N}$ are equivalent. When endowed with the sup norm
and the coordinate wise ordering, $\mathbb{R}^{N}$ can be identified
(algebraically, isometrically and in order) with the Banach lattice $C\left(
\left\{  1,...,N\right\}  \right)  $, where $\left\{  1,...,N\right\}  $
carries the discrete topology.

Suppose that $X$ and $Y$ are two metric spaces and $E$ and $F$ are
respectively ordered vector subspaces (or the positive cones) of
$\mathcal{F}(X)$ and $\mathcal{F}(Y)$ that contain the unity. An operator
$T:E\rightarrow F$ is said to be a \emph{weakly nonlinear operator
}(respectively a\emph{ weakly nonlinear functional }when\emph{ }$F=\mathbb{R}%
$) if it satisfies the following three conditions:

\begin{enumerate}
\item[(SL)] (\emph{Sublinearity}) $T$ is subadditive and positively
homogeneous, that is,%
\[
T(f+g)\leq T(f)+T(g)\quad\text{and}\quad T(af)=aT(f)
\]
for all $f,g$ in $E$ and $a\geq0;$

\item[(M)] (\emph{Monotonicity}) $f\leq g$ in $E$ implies $T(f)\leq T(g).$

\item[(TR)] (\emph{Translatability}) $T(f+\alpha\cdot1)=T(f)+\alpha T(1)$ for
all functions $f\in E$ and all numbers $a\geq0.$
\end{enumerate}

A stronger condition than translatability is that of \emph{comonotonic
additivity},

\begin{enumerate}
\item[(CA)] $T(f+g)=T(f)+T(g)$ whenever the functions $f,g\in E$ are
comonotone in the sense that%
\[
(f(s)-f(t))\cdot(g(s)-g(t))\geq0\text{\quad for all }s,t\in X.
\]

\end{enumerate}

This condition occurs naturally in the context of Choquet's integral (and thus
in the case of Choquet type operators). See \cite{Gal-Nic-2} and
\cite{Gal-Nic-3} and the references therein. For the convenience of the
reader, the basic facts on Choquet's integral are summarized in the Appendix
at the end of this paper.

Of a special interest are the \emph{unital} operators, that is, the\emph{
}operators preserving the unity. A simple example of unital weakly nonlinear
operator is
\[
T:\ell^{\infty}\rightarrow\ell^{\infty},\quad T((x_{n})_{n})=\left(
\limsup\limits_{n\rightarrow\infty}x_{n}\right)  \cdot1;
\]
here $\ell^{\infty}$ is the Banach lattice of all bounded real sequences and
$1$ denotes the sequence with all components equal to unity. As is well known,
$\ell^{\infty}$ can be identified with the space $C_{b}(\mathbb{N})$ (where
$\mathbb{N}$ is endowed with the discrete topology) or with the space
$C(\beta\mathbb{N})$ (where $\beta\mathbb{N}$ is the Stone-Cech
compactification of $\mathbb{N}$). See \cite{Day}. The operator $T$ is not a
Choquet integral (associated to a lower continuous capacity). Indeed,
according to Remark \ref{remApp} $(c)$ in the Appendix,\ $T$ would play the
property%
\[
\lim_{n\rightarrow\infty}T(\xi_{n})=T(\xi),
\]
whenever $(\xi_{n})_{n}$ is a nondecreasing sequence of elements of
$\ell^{\infty}$ that converges coordinatewise to $\xi,$ but this is clearly false.

The permanence properties of weakly nonlinear operators as well as more
examples of such operators are presented in Section 2.

In this paper we extend Korovkin's theorem to the case when the operators
$T_{n}$ are weakly nonlinear operators acting on a function space. This can be
$C(X),$ where $X$ is a compact metric space, or $C_{b}(X),$ where $X$ is a
locally compact metric space. The families of test functions are constructed
via the separating functions (the functions $\gamma:X\times X\rightarrow
\mathbb{R}$ which are continuous and nonnegative and have the property that
$\gamma(x,y)\neq0$ if $x\neq y$). The details are presented in Section 3. This
section also includes the extension of Korovkin's theorem to the case of
weakly nonlinear operators acting on spaces $C(X),$ where $X$ is a compact
space. See Theorem 2. An important consequence of it the following result that
extends our Korovkin type theorem from \emph{\cite{Gal-Nic-1}} in the
particular case of compact metric spaces.

\begin{theorem}
\label{thm1} \emph{(The nonlinear extension of Korovkin's theorem for several
variables)} Suppose that $X$ is a compact subset of the Euclidean space
$\mathbb{R}^{N}$ and let $(T_{n})_{n}$ be a sequence of \textit{sublinear and
monotone operators} from $C(X)$ into itself such that
\begin{equation}
T_{n}(f)(x)\rightarrow f(x)\text{\quad uniformly on }X \label{test}%
\end{equation}
for each of the test functions $1,~\pm\operatorname*{pr}_{1},...,~\pm
\operatorname*{pr}_{N}$ and $\sum_{k=1}^{N}\operatorname*{pr}_{k}^{2}$. Then
\begin{equation}
\lim_{n\rightarrow\infty}T_{n}(f)=f\text{\quad uniformly on }X \label{genconv}%
\end{equation}
for all nonnegative functions $f\in C(X).$ The conclusion \emph{(\ref{genconv}%
)} occurs for all functions $f\in C(X)$ when the operators $T_{n}$ are also translatable.
\end{theorem}

The family of test functions used here is built via the canonical projections
on the Euclidean $N$-dimensional space:
\[
\operatorname*{pr}\nolimits_{k}:(x_{1},...,x_{N})\rightarrow x_{k}%
,\text{\quad}k=1,...,N.
\]

Section 4 is devoted to an extension of Altomare's Korovkin type theorem (see
\cite{Alt2010}, Theorem 3.5, p. 100) to the framework of weakly nonlinear
operators acting on a space $C_{b}(X),$ where $X$ is a locally compact metric space.

Applications of our new results are presented in Section 5.

\section{Preliminaries on weakly nonlinear operators}

Suppose that $X$ and $Y$ are two locally compact compact spaces and $E$ and
$F$ are closed vector sublattices respectively of the Banach lattices
$C_{b}(X)$ and $C_{b}(Y).$

Every monotone and subadditive operator $T$ $:E\rightarrow F$ verifies the
inequality%
\begin{equation}
\left\vert T(f)-T(g)\right\vert \leq T\left(  \left\vert f-g\right\vert
\right)  \text{\quad for all }f,g. \label{f1}%
\end{equation}
Indeed, $f\leq g+\left\vert f-g\right\vert ~$yields $T(f)\leq T(g)+T\left(
\left\vert f-g\right\vert \right)  ,$ that is, $T(f)-T(g)\leq T\left(
\left\vert f-g\right\vert \right)  $, and interchanging the role of $f$ and
$g$ we infer that $-\left(  T(f)-T(g)\right)  \leq T\left(  \left\vert
f-g\right\vert \right)  .$

If $T$ is linear, then the property of monotonicity is equivalent to that of
\emph{positivity}, that is, to the fact that%
\[
T(f)\geq0\text{\quad for all }f\geq0.
\]
If the operator $T$ is monotone and positively homogeneous, then necessarily%
\[
T(0)=0.
\]
Every sublinear operator is convex and a convex function $\Phi:E\rightarrow F$
is sublinear if and only if it is positively homogeneous.

The continuity of a sublinear operator $T:E\rightarrow F$ is equivalent to its
continuity at the origin, which in turn is equivalent to existence of a
constant $\lambda\geq0$ such that%
\[
\left\Vert T\left(  x\right)  \right\Vert \leq\lambda\left\Vert x\right\Vert
\text{\quad for all }x\in E.
\]
The smallest constant $\lambda=\left\Vert T\right\Vert $ with this property
will be called the \emph{norm} of $T.$

\begin{remark}
If $T:C_{b}(X)\rightarrow C_{b}(X)$ is a sublinear and monotone operator, then
$T$ is continuous and%
\[
\left\Vert T\right\Vert =\left\Vert T(1)\right\Vert .
\]
Indeed, $\left\vert f\right\vert \leq\left\Vert f\right\Vert _{\infty}\cdot1,$
so that, according to \emph{(\ref{f1})}, we infer that
\[
\left\Vert T(f)\right\Vert \leq\left\Vert f\right\Vert _{\infty}\left\Vert
T(1)\right\Vert .
\]
This shows that $T$ is continuous and $\left\Vert T\right\Vert \leq\left\Vert
T(1)\right\Vert ;$ the other inequality is trivial an thus $\left\Vert
T\right\Vert =\left\Vert T(1)\right\Vert .$ An immediate consequence is that
$\left\Vert T\right\Vert =1$ when $T$ is in addition unital.
\end{remark}

The following variant of H\"{o}lder's inequality is a particular case of
Theorem 3 in our paper \cite{Gal-Nic-2}.

\begin{lemma}
\label{lemHolder}$($\emph{H\"{o}lder's inequality for} $p\in(1,\infty)$ and
$1/p+1/q=1)$ Suppose that $X$ is a compact metric space\textit{ and
}$T:C(X)\rightarrow C(X)$\textit{ }is a unital, \textit{sublinear and monotone
operator}.\textit{ Then}%
\[
T(|fg|)\leq\lbrack T(|f|^{p})]^{1/p}\cdot\lbrack T(|g|^{q})]^{1/q}.
\]
for all $f,g\in C(X).$
\end{lemma}

Concrete examples of sublinear and monotone operators are presented in
\cite{NO2020} and references therein. They are ubiquitous in many fields like
functional analysis, convex analysis and partial differential equations. As we
prove in Section 3, even the sequences of sublinear and monotone operators
$T_{n}:C(X)\rightarrow C(X)$ having the property that $T_{n}(1)=1$ for every
$n\in\mathbb{N}$ offer a natural framework for approximating the nonnegative
continuous functions by suitable special classes of functions.

The set $\mathcal{WN}(C_{b}(X),C_{b}(X))$ of all weakly nonlinear operators
$T:C_{b}(X)\rightarrow C_{b}(X)$ is a convex cone in the Banach space
$\operatorname*{Lip}_{0}\left(  C_{b}(X),C_{b}(X)\right)  ,$ of \ all
Lipschitz maps from $C_{b}(X)$ into itself that vanish at the origin. In turn,
it includes the cone $\mathcal{L}_{+}\left(  C_{b}(X),C_{b}(X)\right)  ,$ of
all linear, continuous and monotone operators from $C_{b}(X)$ into itself.

A general procedure to generate new weakly nonlinear operators form old ones
is as follows:

\begin{lemma}
\label{lemsup}$(a)$ If $S,T\in\mathcal{WN}(C_{b}(X),C_{b}(X))$ and
$S(1)=T(1),$ then the operator $S\vee T$ defined by the formula
\[
\left(  S\vee T\right)  (f)=\sup\left\{  S(f),T(f)\right\}  \quad\text{for
}f\in C(X),
\]
also belongs to $\mathcal{WN}(C(X), C(X)).$

$(b)$ If $S,T\in\mathcal{WN}(C_{b}(X),C_{b}(X))$ and $T$ is unital, then
$ST\in\mathcal{WN}(C_{b}(X),C_{b}(X)).$
\end{lemma}

\begin{proof}
Indeed, the fact that the pointwise sup of two sublinear and monotone
operators is also sublinear and monotone is obvious. In addition,%
\begin{align*}
\left(  S\vee T\right)  (f+\alpha1)  &  =\sup\left\{  S(f+\alpha
1),T(f+\alpha1)\right\} \\
&  =\sup\left\{  S(f)+\alpha S(1),T(f)+\alpha T(1)\right\} \\
&  =\sup\left\{  S(f),T(f)\right\}  + \sup\{\alpha S(1), \alpha T(1)\}\\
&  =\sup\left\{  S(f),T(f)\right\}  +\alpha\left(  S\vee T\right)  (1)
\end{align*}
for all $f\in C(X)$ and $\alpha\ge0.$ The proof is done.
\end{proof}

\begin{example}
\label{ex1}The Banach lattice $c$ of all convergent sequences of real numbers
\emph{(}endowed with the sup norm and the coordinatewise ordering\emph{)} can
be identified with $C(\mathbb{\hat{N}}),$ where $\mathbb{\hat{N}=N}%
\cup\{\infty\}$ is the one point compactification of the discrete space
$\mathbb{N}.$ See \emph{\cite{Day}}. According to Lemma \emph{\ref{lemsup}}
$(a),$ the following operators, from $c$ into itself, are unital and weakly
nonlinear:%
\begin{align*}
T_{1}\left(  (x_{n})_{n}\right)   &  =\left(  \sup\{x_{n},\lim_{k\rightarrow
\infty}x_{k}\}\right)  _{n}\\
T_{2}\left(  (x_{n})_{n}\right)   &  =\left(  \sup\{\frac{x_{1}+\cdots+x_{n}%
}{n},\lim_{k\rightarrow\infty}x_{k}\}\right)  _{n}\\
T_{3}\left(  (x_{n})_{n}\right)   &  =\left(  \sup\{x_{n},\frac
{x_{1}+2x_{2}+\cdots+2^{n}x_{n}}{2^{n+1}-1}\}\right)  _{n}.
\end{align*}

\end{example}

\section{The case of compact metric spaces}

The basic ingredient in our approach of extending Korovkin's theory is a
technical estimate for uniformly continuous functions, originating in his
paper \cite{Ko1953} from 1953, and put here in a slightly more generality.

\begin{lemma}
\label{lemgamma}If $X=(X,d)$ is a compact metric space, and $\gamma:X\times
X\rightarrow\mathbb{R}$ is a separating function, that is, a nonnegative
continuous function such that%
\[
\gamma(x,y)=0\text{\quad implies }x=y,
\]
then every real-valued continuous function $f$ defined on $X$ verifies an
estimate of the form%
\[
\left\vert f(x)-f(y)\right\vert \leq\varepsilon+\delta(\varepsilon
)\gamma(x,y)\text{\quad for all }x,y\in X\text{ and }\varepsilon>0.
\]

\end{lemma}

\begin{proof}
We borrow the quick argument from \cite{N1979}, \cite{N2009}. If the estimate
above doesn't work, then for a suitable $\varepsilon_{0}>0$ one can find two
sequences $(x_{n})_{n}$ and $(y_{n})_{n}$ of elements of $X$ such that%
\begin{equation}
\left\vert f(x_{n})-f(y_{n})\right\vert \geq\varepsilon_{0}+2^{n}\gamma
(x_{n},y_{n}) \label{st}%
\end{equation}
for all $n.$ Without loss of generality we may assume (by passing to
subsequences) that both sequences $(x_{n})_{n}$ and $(y_{n})_{n}$ are
convergent, respectively to $x$ and $y.$ Since $f$ is bounded, the inequality
(\ref{st}) forces $x=y.$ Indeed,%
\[
\frac{\left\vert f(x_{n})-f(y_{n})\right\vert }{2^{n}}\rightarrow0\text{ and
}\frac{\left\vert f(x_{n})-f(y_{n})\right\vert }{2^{n}}\geq\gamma(x_{n}%
,y_{n})\rightarrow\gamma(x,y)\geq0,
\]
which implies that $\gamma(x,y)=0.$ On the other hand, from (\ref{st}) one can
infer that $\left\vert f(x)-f(y)\right\vert \geq\varepsilon_{0}$ and thus
$x\neq y.$ This contradiction shows that the assumption made at the beginning
of the proof is wrong and the assertion of Lemma \ref{lemgamma} is true.
\end{proof}

\begin{remark}
\label{rem0}The argument of Lemma \emph{\ref{lemgamma}} also shows that every
separating function $\gamma:X\times X\rightarrow\mathbb{R}$ is related to the
metric $d$ on $X$ via an estimate of the form%
\[
d(x,y)\leq\varepsilon+\delta(\varepsilon)\gamma(x,y)\text{\quad for all
}(x,y)\in X\times X\text{ and }\varepsilon>0.
\]

\end{remark}

If $X=(X,d)$ is an arbitrary compact metric space and $\varphi:[0,\infty
)\rightarrow\lbrack0,\infty)$ is a continuous function such that
$\varphi(t)>0$ for $t>0,$ then $\gamma(x,y)=\varphi(d(x,y))$ is an example of
separating function.

Every separating function $\gamma:X\times X\rightarrow\mathbb{R}$ generates a
family of nonnegative continuous functions on $X,$ precisely,%
\[
\gamma_{x}:X\rightarrow\mathbb{R},\text{\quad}\gamma_{x}(y)=\gamma
(x,y)~\text{for }x,y\in X.
\]
As shows Theorem \ref{thm2} below, this family can be used as a family of test
functions in the same manner as the functions $1,~x$ and $x^{2}$ were used in
Korovkin's theorem. Therefore we are primarily interested in separating
functions producing a minimal number of test functions. A classical example is
offered by the case of compact subsets $X$ of $\mathbb{R}^{N}.$ Choosing
$f_{1},...,f_{m}\in C(X)$ a family of functions which separates the points of
$X$ and
\begin{equation}
\gamma(x,y)=\sum_{k=1}^{m}\left(  f_{k}(x)-f_{k}(y)\right)  ^{2} \label{csf}%
\end{equation}
is a separating function; when this family consists of the coordinate
functions $\operatorname*{pr}_{1},...,\operatorname*{pr}_{N},$ then%
\[
\gamma(x,y)=\left\Vert x-y\right\Vert ^{2}.
\]
The following result represents a nonlinear generalization of Korovkin's
theorem and of many other related results existing in the literature.

\begin{theorem}
\label{thm2}Let $X$ be a compact metric space $($endowed with the metric $d)$
and let $(T_{n})_{n}$ be a sequence of sublinear and monotone operators from
$C(X)$ into itself such that
\begin{equation}
T_{n}(1)(x)\rightarrow1\text{\quad uniformly on }X. \label{hyp1}%
\end{equation}
Suppose that $\gamma:X\times X\rightarrow\mathbb{R}$ is a separating function
such that
\begin{equation}
T_{n}(\gamma_{x})(x)\rightarrow0\text{\quad uniformly on }X. \label{hyp2}%
\end{equation}
Then for all nonnegative functions $f\in C(X),$%
\begin{equation}
T_{n}(f)\rightarrow f~\text{\quad uniformly on }X. \label{conclusion}%
\end{equation}
This convergence occurs for all $f\in C(X)$ if the operators \thinspace$T_{n}$
are also translatable $($that is, when they are weakly nonlinear$)$.
\end{theorem}

\begin{proof}
Let $f\in C(X)$ be a nonnegative function. Then, according to Lemma
\ref{lemgamma}, for every $\varepsilon>0$ there is $\delta(\varepsilon)>0$
such that%
\begin{equation}
\left\vert f-f(x)\right\vert \leq\varepsilon+\delta(\varepsilon)\gamma_{x}
\label{absgamma}%
\end{equation}
\quad for all $x\in X.$ Since the operators $T_{n}$ are subadditive,
positively homogeneous and monotone, one can use the inequality (\ref{f1}) to
show that%
\begin{align*}
\left\vert T_{n}(f)-f(x)T_{n}(1)\right\vert  &  =\left\vert T_{n}%
(f)-T_{n}(f(x)\cdot1)\right\vert \leq T_{n}\left(  \left\vert
f-f(x)\right\vert \right) \\
&  \leq\varepsilon T_{n}(1)+\delta(\varepsilon)T_{n}(\gamma_{x}).
\end{align*}
According to our hypotheses (\ref{hyp1}) and (\ref{hyp2}), this leads to the
conclusion that $T_{n}(f)\rightarrow f~,$ uniformly on $X.$

Suppose now that each operator $T_{n}$ is also weakly additive. Every function
$f\in C(X)$ verifies the inequality $f(x)+\left\Vert f\right\Vert _{\infty
}\geq0$, whenever $x\in X$, so that by taking into account the above
considerations, we infer that
\[
T_{n}(f+\left\Vert f\right\Vert _{\infty})(x)\rightarrow f(x)+\left\Vert
f\right\Vert _{\infty},
\]
uniformly on $X$. Taking into account the hypothesis (\ref{hyp1}) and the fact
that the operators $T_{n}$ were assumed to be weakly additive, we have
\[
T_{n}(f+\left\Vert f\right\Vert _{\infty})(x)=T_{n}(f)(x)+\left\Vert
f\right\Vert _{\infty}\cdot T_{n}(1)(x)\rightarrow T_{n}(f)(x)+\left\Vert
f\right\Vert _{\infty},
\]
which yields that $T_{n}(f)\rightarrow f~,$ uniformly on $X,$ for any function
$f\in C(X).$ The proof is done.
\end{proof}

According to Remark \ref{rem0}, if the sequence of operators $T_{n}$ verifies
the condition (\ref{hyp2}) in Theorem \ref{thm2}, it also verifies the
condition%
\[
T_{n}(d_{x})(x)\rightarrow0\text{\quad uniformly on }X.
\]
This outlines the prominent role played by the distance function among the
separating functions.

\begin{theorem}
\label{thm3}Under the hypotheses of Theorem \emph{\ref{thm2}}, if $\gamma=d$
and $f\in C(X)$ is a Lipschitz continuous function with the Lipschitz constant
$K,$ then the following estimate holds:
\[
|T_{n}(f)(x)-f(x)|\leq K\cdot\sup\left\{  \left\vert T_{n}(d_{x}%
^{2})(x)\right\vert ^{1/2}:x\in X\right\}  \text{\quad for all }x\in X\text{
and }n\in\mathbb{N}.
\]

\end{theorem}

\begin{proof}
The argument is similar to that of Theorem \ref{thm2}, replacing the starting
estimate (\ref{absgamma}) by the condition of Lipschitzianity,
\[
f(x)\cdot1-K\cdot d_{x}\leq f\leq f(x)\cdot1+K\cdot d_{x}\text{\quad for all
}x\in X.
\]
Suppose for a moment that $f\geq0.$ Since the operators $T_{n}$ are
subadditive, monotonic and positively homogeneous one can apply them to the
left-hand side inequality (rewritten as $f(x)\cdot1\leq f+Kd_{x}),$ resulting
that%
\[
f(x)\leq T_{n}(f)(x)+KT_{n}(d_{x})(x).
\]
Applying these operators to the right hand side inequality one obtains
\[
T_{n}(f)(x)\leq f(x)+KT_{n}(d_{x})(x).
\]
Therefore, taking into account Lemma \ref{lemHolder}, we conclude that
\[
|T_{n}(f)(x)-f(x)|\leq KT_{n}(d_{x})(x)\leq K\left(  T_{n}(d_{x}%
^{2})(x)\right)  ^{1/2}\text{\quad for all }x\in X\text{ and }n\in\mathbb{N}.
\]

The case of Lipschitz functions not necessarily nonnegative can be settled as
in the proof of Theorem \ref{thm2}.
\end{proof}

\begin{proof}
[Proof of Theorem 1.]When $X$ is a compact subset of $\mathbb{R}^{N}$ and
$\gamma(x,y)=\left\Vert x-y\right\Vert ^{2},$ Theorem \emph{\ref{thm2}} can be
restated in a more convenient way by replacing the two conditions
\emph{(\ref{hyp1})}$\&$\emph{(\ref{hyp2})} with a set of $2N+1$ tests of
convergence:
\[
T_{n}(f)(x)\rightarrow1\text{\quad uniformly on }X,
\]
for each of the test functions $1,~\pm\operatorname*{pr}\nolimits_{1}%
,...,~\pm\operatorname*{pr}\nolimits_{N}~$and $\sum_{k=1}^{N}%
\operatorname*{pr}\nolimits_{k}^{2}.$ Here we can replace $\sum_{k=1}%
^{N}\operatorname*{pr}\nolimits_{k}^{2}$ by the string of test functions
$\operatorname*{pr}\nolimits_{1}^{2},...,\operatorname*{pr}\nolimits_{N}^{2}$.

Indeed, by denoting%
\[
M=\sup_{x\in X}\left\{  \operatorname*{pr}\nolimits_{1}%
(x),...,\operatorname*{pr}\nolimits_{N}(x),0\right\}  ,
\]
we have
\begin{multline*}
0\leq T_{n}(\Vert\cdot-x\Vert^{2})(x)\leq T_{n}(\left\Vert x\right\Vert
^{2})(x)+2T_{n}(-\langle\cdot,x\rangle)(x)+\left\Vert x\right\Vert ^{2}%
T_{n}(1)(x)\\
=T_{n}\left(  \sum_{k=1}^{N}\operatorname*{pr}\nolimits_{k}^{2}\right)
(x)+2T_{n}\left[  \sum_{k=1}^{N}\left(  -\operatorname*{pr}\nolimits_{k}%
(x)\right)  \cdot\operatorname*{pr}\nolimits_{k}(\cdot)\right]  (x)+\left\Vert
x\right\Vert ^{2}T_{n}(1)(x)\\
=T_{n}\left(  \sum_{k=1}^{N}\operatorname*{pr}\nolimits_{k}^{2}\right)
(x)+2T_{n}\left[  \sum_{k=1}^{N}\left(  M-\operatorname*{pr}\nolimits_{k}%
(x)\right)  \cdot(\operatorname*{pr}\nolimits_{k}(\cdot))+M\sum_{k=1}%
^{N}\left(  -\operatorname*{pr}\nolimits_{k}(\cdot)\right)  \right]  (x)\\
+\left\Vert x\right\Vert ^{2}T_{n}(1)(x)\\
\leq T_{n}\left(  \sum_{k=1}^{N}\operatorname*{pr}\nolimits_{k}^{2}\right)
(x)+2\sum_{k=1}^{N}\left(  M-\operatorname*{pr}\nolimits_{k}(x)\right)  \cdot
T_{n}\left(  \operatorname*{pr}\nolimits_{k}(\cdot)\right)  (x)+2M\sum
_{k=1}^{N}T_{k}\left(  -\operatorname*{pr}\nolimits_{k}(\cdot)\right) \\
+\left\Vert x\right\Vert ^{2}T_{n}(1)(x)
\end{multline*}
and assuming that $\lim_{n\rightarrow\infty}T_{n}(f)(x)\rightarrow f$\quad
uniformly on $X$ for each of the test functions $1,~\pm\operatorname*{pr}%
\nolimits_{1},...,~\pm\operatorname*{pr}\nolimits_{N}~$and $\sum_{k=1}%
^{N}\operatorname*{pr}\nolimits_{k}^{2}$ we can easily check that the
right-hand side of the precedent string of inequalities converges uniformly to
$0$ on $X.$ Consequently $T_{n}(\Vert\cdot-x\Vert^{2})(x)\rightarrow0$
uniformly on $X$ and Theorem \emph{\ref{thm2}} applies. The proof is done.
\end{proof}

\begin{remark}
Working with a finite family $f_{1},...,f_{p}$ of continuous functions that
separates the points of $X$ and the separating function $\gamma(x,y)=\sum
_{k=1}^{p}\left\vert f_{k}(x)-f_{k}(y)\right\vert ^{2}$, one can arrive at the
conclusion of Theorem 1 by verifying the convergence \emph{(\ref{conclusion})}
for the $\ 2p+2$ test functions, $1,~\pm f_{1},...,~\pm f_{p}~$and $\sum
_{k=1}^{2}f_{k}^{2}.$
\end{remark}

\begin{example}
Consider now the particular case of the unit circle
\[
S^{1}=\left\{  \left(  \cos\varphi,\sin\varphi\right)  :\varphi\in
\mathbb{R}\right\}  .
\]
With respect to the metric induced by $\mathbb{R}^{2},$%
\begin{align*}
d(\left(  \cos\varphi,\sin\varphi\right)  ,\left(  \cos\psi,\sin\psi\right)
)  &  =\sqrt{\left(  \cos\varphi-\cos\psi\right)  ^{2}+\left(  \sin
\varphi-\sin\psi\right)  ^{2}}\\
&  =2\left\vert \sin\frac{\varphi-\psi}{2}\right\vert ,
\end{align*}
$S^{1}$ is a compact $($that is, bounded and close$)$ subset of $\mathbb{R}%
^{2}.$ Choosing as a separating function the square distance,
\begin{align*}
\gamma\left(  \left(  \cos\varphi,\sin\varphi\right)  ,\left(  \cos\psi
,\sin\psi\right)  \right)   &  =\sin^{2}\frac{\varphi-\psi}{2}\\
&  =1-\cos\varphi\cos\psi-\sin\varphi\sin\psi
\end{align*}
one can easily check that the conditions \emph{(\ref{hyp1})}$\&$%
\emph{(\ref{hyp2}) }in\emph{ }Theorem\emph{ }$2$ can be replaced in this case
by the fulfillment of the following $5$ tests of convergence:
\[
T_{n}(f)(x)\rightarrow1\text{\quad uniformly on }X
\]
for each of the functions $1,~\pm\operatorname*{pr}\nolimits_{1}$ and
$\pm\operatorname*{pr}\nolimits_{2}.$ It is well known that the Banach space
$C(S^{1})$ $,$ can be identified with the space $C_{2\pi}(\mathbb{R})$, of all
continuous and $2\pi$-periodic functions $f:\mathbb{R\rightarrow R}.$ Modulo
this identification, we infer from Theorem $1$ that a sufficient condition for
a sequence of weakly nonlinear operators $T_{n}:C_{2\pi}(\mathbb{R}%
)\rightarrow C_{2\pi}(\mathbb{R})$ to verify the condition
\[
\lim_{n\rightarrow\infty}T_{n}(f)(\varphi)\rightarrow0,\text{\quad uniformly
on }\mathbb{R}%
\]
is to verify this conditions for the test functions $1,~\pm\cos\varphi~$and
$\pm\sin\varphi.$ This was first noticed by Korovkin \emph{\cite{Ko1953}},
\emph{\cite{Ko1960}} in the particular case of linear operators.
\end{example}

\begin{example}
\label{exPopa}Suppose that $X=(X,d)$ is a compact metric space. Then the
product space $X\times S^{1}$ is also a compact metric space and the space
$C(X\times S^{1})$ can be identified with the Banach space $C_{2\pi}%
(K\times\mathbb{R}),$ of all continuous functions $f:K\times\mathbb{R}%
\rightarrow\mathbb{R}$, $2\pi$-periodic in the second variable, endowed with
the sup norm. This space is genuine for many results in dynamical systems
theory. By considering the separating function
\[
\gamma((x,\varphi),(y,\psi))=d(x,y)^{2}+\sin^{2}\frac{\varphi-\psi}{2},
\]
Popa \emph{\cite{Popa}} has recently proved the variant of Theorem
\emph{\ref{thm1} }for the linear and positive operators\emph{ }$T:C(X\times
S^{1})\rightarrow C(X\times S^{1}).$ The reader can easily check that actually
his results extend to the case of weakly nonlinear operators. In particular,
when $X$ is compact subset of $\mathbb{R}^{N}$ then the convergence%
\[
T_{n}(f)\rightarrow f~\text{\quad uniformly on }X\times S^{1},
\]
for all $f\in C(X\times S^{1})$ reduces to its verification for the product
functions $f(x)=u(x)v(\varphi),$ where
\[
u\in\left\{  1,~\pm\operatorname*{pr}\nolimits_{1},...,~\pm\operatorname*{pr}%
\nolimits_{N}~\text{and}\sum_{k=1}^{N}\operatorname*{pr}\nolimits_{k}%
^{2}\right\}  \text{ and }v\in\left\{  1,~\pm\cos\varphi~,~\pm\sin
\varphi\right\}  .
\]

We left to the reader the easy exercise to detail Theorem $1$ in some other
cases of interest such as the torus $S^{1}\times S^{1}$ and the $2$%
-dimensional sphere $S^{2}$.
\end{example}

\begin{remark}
In the absence of the condition of translatability the conclusion of Theorem
\emph{\ref{thm2}} may fail for functions with variable sign. An example
working for $X=[0,1]$ is given by the Bernstein like operators
\[
T_{n}(f)(x)=\sum_{k=0}^{n}\binom{n}{k}x^{k}(1-x)^{n-k}\sup\left\{  f\left(
k/n\right)  ,0\right\}  ,
\]
which are sublinear and monotone (but not translatable). Clearly,
$T_{n}(f)\rightarrow f$ uniformly on $[0,1]$ for each of the functions $1$,
$x$ and $x^{2}.$ According to Theorem \emph{\ref{thm2}, }this convergence
occurs for all nonnegative functions $f\in C\left(  [0,1]\right)  .$ Clearly,
it fails for the nonpositive functions. Remarkably, the case of nonnegative
functions is strong enough to provide valuable information for \emph{all
}functions in\emph{ }$C\left(  [0,1]\right)  ,$ for example, the possibility
to approximate them by polynomials. Indeed,%
\[
T_{n}(f+\left\Vert f\right\Vert _{\infty})-\left\Vert f\right\Vert _{\infty
}\rightarrow f\text{\quad uniformly on }[0,1].
\]

\end{remark}

\section{The extension of a result due to Altomare}

The next theorem represents a nonlinear analogue of a result due to Altomare
(see \cite{Alt2010}, Theorem 3.5, p. 100).

\begin{theorem}
\label{thm4}Let $X$ be a locally compact metric space $($endowed with the
metric $d)$ and consider a vector sublattice $E$ of $\mathcal{F}(X)$
containing the constant functions and all the functions $d_{x}^{p}$ for $x\in
X$ and some exponent $p\geq1$. Let $(T_{n})_{n}$ be a sequence of sublinear
and monotone operators from $E$ into $\mathcal{F}(X)$ which verifies the
following two conditions:

$(a)$ $\lim_{n\rightarrow\infty}T_{n}(1)=1$, uniformly on compact subsets of
$X;$

$(b)$ $\lim_{n\rightarrow\infty}T_{n}(d_{x}^{p})(x)=0$, uniformly on compact
subsets of $X;$

Then, for all nonnegative $f$ in $E\cap C_{b}(X)$, we have%
\[
\lim_{n\rightarrow\infty}T_{n}(f)=f,\text{\quad uniformly on compact subsets
of }X.
\]
The convergence occurs for all functions in $E\cap C_{b}(X)$ when the
operators $T_{n}$ are also translatable.
\end{theorem}

The proof of Theorem 4, needs the following lemma due to Altomare. See
\cite{Alt2010}, Lemma 3.4, p. 99 for details.

\begin{lemma}
\label{lemAlt}Let $X$ be a locally compact metric space endowed with the
metric $d$. Then for every compact subset $K$ of $X$ and for every
$\varepsilon>0$, there exist $0<\varepsilon^{\prime}<\varepsilon$ and a
compact subset $K_{\varepsilon}$ of $X$ such that the open ball
$B_{\varepsilon^{\prime}}(x)$ is included in $K_{\varepsilon}$ for every $x\in
K$.
\end{lemma}

\begin{proof}
[Proof of Theorem $\mathbf{4}$.]Let $f\in E\cap C_{b}(X)$ and $\varepsilon>0$
be arbitrarily fixed. Then for every compact subset $K$ of $X,$ choose
$\varepsilon^{\prime}\in(0,\varepsilon)$ and $K_{\varepsilon}$ be as in Lemma
\emph{\ref{lemAlt}}.

Since $f$ is uniformly continuous on $K_{\varepsilon}$, there exists
$\delta\in(0,\varepsilon^{\prime})$ such that
\[
|f(x)-f(y)|\leq\varepsilon\text{\quad for every }x,y\in K_{\varepsilon}\text{
with }d(x,y)\leq\delta.
\]
Suppose that $x\in K$ and $y\in X$. If $d(x,y)\leq\delta$, then $y\in
B^{\prime}(x,\overline{\varepsilon})\subset K_{\varepsilon}$ and therefore,
$|f(x)-f(y)|\leq\varepsilon$. If $d(x,y)\geq\delta$, then
\[
|f(x)-f(y)|\leq\frac{2\Vert f\Vert_{\infty}}{\delta^{p}}\cdot d^{p}(x,y).
\]
Therefore
\[
|f-f(x)|\leq\frac{2\Vert f\Vert_{\infty}}{\delta^{p}}\cdot d_{x}%
^{p}+\varepsilon\cdot1\text{ for all }x\in K.
\]
so that, taking into account the inequality (\ref{f1}), we infer in the case
of nonnegative functions $f$ that
\begin{equation}
|T_{n}(f)(x)-f(x)T_{n}(1)(x)|=T_{n}(|f-f(x)|)(x)\leq\frac{2\Vert
f\Vert_{\infty}}{\delta^{p}}\cdot d_{x}^{p}(x)+\varepsilon T_{n}(1)(x)
\label{estth2}%
\end{equation}
for all $x\in K,$ whence%
\begin{equation}
\lim_{n\rightarrow\infty}T_{n}(f)(x)=f(x),\text{ uniformly with respect to
}x\in K. \label{estth2+}%
\end{equation}
Assume now that all operators $T_{n}$ are translatable. Replacing $f$\ by
$f+\Vert f\Vert_{\infty},$ the equality in the left hand side (\ref{estth2})
still works, so that the inequality in the right hand side occurs for all
functions $f\in E\cap C_{b}(X).$ The same is true concerning the formula
(\ref{estth2+}) and the proof is done.
\end{proof}

An example illustrating Theorem 4 is exhibited at the end of the next section.

\section{Applications}

In this section we illustrate the results in the previous sections by several
concrete examples. We adopt the convention $0^{0}=1$.

\medskip

\noindent\textbf{The Bernstein-Kantorovich-Choquet polynomial operators}. We
proved in \cite{Gal-Nic-1} that the Bernstein-Kantorovich-Choquet polynomial
operators for functions of one real variable,%
\[
K_{n,\mu}^{(1)}:C([0,1])\rightarrow C([0,1]),
\]
defined by the formula
\[
K_{n,\mu}^{(1)}(f)(x)=\sum_{k=0}^{n}p_{n,k}(x)\cdot\frac{(C)\int
_{k/(n+1)}^{(k+1)/(n+1)}f(t)\mathrm{d}\mu(t)}{\mu([k/(n+1),(k+1)/(n+1)])},
\]
verifies the conditions $K_{n,\mu}^{(1)}(x^{k})\rightarrow x^{k}$ uniformly on
$[0,1]$ for $k\in\{0,1,2\},$ which implies that $K_{n,\mu}^{(1)}(f)\rightarrow
f$ uniformly on $[0,1]$ for all functions $f\in C([0,1]).$

The Bernstein-Kantorovich-Choquet polynomial operators for functions of two
real variables,%
\[
K_{n,\mu}^{(2)}:C([0,1]^{2})\rightarrow C([0,1]^{2}),
\]
are defined by the formula
\begin{multline*}
K_{n,\mu}^{(2)}(f)(x_{1},x_{2})=\sum_{k_{1}=0}^{n}\sum_{k_{2}=0}^{n}%
p_{n,k_{1}}(x_{1})p_{n,k_{2}}(x_{2})\\
\cdot\frac{(C)\int_{k_{1}/(n+1)}^{(k_{1}+1)/(n+1)}\left(  (C)\int
_{k_{2}/(n+1)}^{(k_{2}+1)/(n+1)}f(t_{1},t_{2})\mathrm{d}\mu(t_{2})\right)
\mathrm{d}\mu(t_{1})}{\mu([k_{1}/(n+1),(k_{1}+1)/(n+1)])\mu([k_{2}%
/(n+1),(k_{2}+1)/(n+1)])},
\end{multline*}
where%
\[
p_{n,k}(t)={\binom{n}{k}}t^{k}(1-t)^{n-k},\text{\quad for }t\in\lbrack
0,1]\text{ and }n\in\mathbb{N},
\]
$\mu=\sqrt{\mathcal{L}}$ is the monotone and submodular $($and therefore
subadditive$)$ set function associated to the Lebesgue measure $\mathcal{L}$
on the specific interval of integration and $f\in C([0,1]^{2})$.

Due to the properties of the Choquet integral mentioned in the Appendix, it
follows that each operator $K_{n,\mu}^{(2)}$ is a weekly nonlinear and unital
operator from $C([0,1]^{2})$ into itself. However, this operator is not
comonotonically additive (and thus escapes the theory developed in
\cite{Gal-Nic-1}).

We will show that%
\begin{equation}
K_{n,\mu}^{(2)}(f)(x_{1},x_{2})\rightarrow f(x_{1},x_{2})\text{\quad uniformly
on }[0,1]^{2} \label{convBKC}%
\end{equation}
for all test functions \thinspace$1,~\pm\operatorname*{pr}\nolimits_{1}%
,~\pm\operatorname*{pr}\nolimits_{2},~\operatorname*{pr}\nolimits_{1}%
^{2}+\operatorname*{pr}\nolimits_{2}^{2}$ $($which will imply, via Theorem
$1$, that this convergence occurs for all functions $f\in C([0,1]^{2}).$

The case of the unity is clear, while the case of the functions $\pm
\operatorname*{pr}\nolimits_{1}$and $\pm\operatorname*{pr}\nolimits_{2}$ is
the settled by the aforementioned properties of the operators $K_{n,\mu}%
^{(1)}.$ As concerns the case of the function $\operatorname*{pr}%
\nolimits_{1}^{2}+\operatorname*{pr}\nolimits_{2}^{2},$ notice that
\begin{multline*}
K_{n,\mu}^{(2)}(|\operatorname*{pr}\nolimits_{1}^{2}(\mathbf{t}%
)+\operatorname*{pr}\nolimits_{2}^{2}(\mathbf{t})-\operatorname*{pr}%
\nolimits_{1}^{2}(\mathbf{x})+\operatorname*{pr}\nolimits_{2}^{2}%
(\mathbf{x})|)\\
\leq\sum_{i=1}^{2}K_{n,\mu}^{(2)}(|\operatorname*{pr}\nolimits_{i}%
(\mathbf{t})+\operatorname*{pr}\nolimits_{i}(\mathbf{x})|\cdot
|\operatorname*{pr}\nolimits_{i}(\mathbf{t})-\operatorname*{pr}\nolimits_{i}%
(\mathbf{x})|)\\
\leq2\sum_{i=1}^{2}K_{n,\mu}^{(2)}(|\operatorname*{pr}\nolimits_{i}%
(\mathbf{t})-\operatorname*{pr}\nolimits_{i}(\mathbf{x})|)\\
\leq2\sum_{i=1}^{2}\sqrt{K_{n,\mu}^{(2)}(|\operatorname*{pr}\nolimits_{i}%
(\mathbf{t})-\operatorname*{pr}\nolimits_{i}(\mathbf{x})|^{2})}\\
=2\sum_{i=1}^{2}\sqrt{K^{(2)}_{n,\mu}(t_{i}^{2})(x)+2x_{i}K^{(2)}_{n,\mu
}(-t_{i})(x)+x_{i}^{2}}%
\end{multline*}
according to Lemma \ref{lemHolder}. Now, by using the calculations for the
Bernstein-Kantorovich-Choquet operators in one variable in \cite{Gal-1},
\cite{Gal-Nic-1}, it is immediate that $K_{n,\mu}^{(2)}(t_{i}^{2}%
)(\mathbf{x})\rightarrow x_{i}^{2}$ and $K_{n,\mu}^{(2)}(-t_{i})(\mathbf{x}%
)\rightarrow-x_{i}$ as $n\rightarrow\infty$, uniformly with respect to
$\mathbf{x}=(x_{1},x_{2})\in\lbrack0,1]^{2}$, for $i\in\{1,2\}.$ Therefore, it
follows that
\[
K^{(2)}_{n,\mu}(t_{i}^{2})(\mathbf{x})+2x_{i}K^{(2)}_{n,\mu}(-t_{i}%
)(\mathbf{x})+x_{i}^{2}\rightarrow0,
\]
uniformly with respect to $\mathbf{x},$ for $i\in\{1,2\}.$ Thus the
convergence (\ref{convBKC}) also occurs for the function $\operatorname*{pr}%
\nolimits_{1}^{2}+\operatorname*{pr}\nolimits_{2}^{2}$ (and thus for all
functions $f\in C([0,1]^{2}).$

The reader can now easily extend this example to the case of
Bernstein-Kantorovich-Choquet polynomial operators for functions of $N$ real variables.

Notice that while $K_{n}^{(2)}$ is only translatable, the one variable
corresponding operator $K_{n}^{(1)}$ is comonotonic additive, see
\cite{Gal-Nic-1}. Also, in the one variable case, the error estimate in
approximation of $f$ by $K^{(1)}_{n}(f)$ in terms of the modulus of continuity
was obtained in \cite{Gal-1}.

\medskip

\noindent\textbf{The bivariate possibilistic Bernstein-Durrmeyer and
Kantorovich polynomial operators.} In the case of one variable, the so-called
possibilistic Bernstein-Durrmeyer polynomials operators and possibilistic
Kantorovich polynomial operators were considered in \cite{Gal-3} by replacing
in the expressions of the classical integral operators of Bernstein-Durrmeyer
and of Kantorovich, the Lebesgue integral by the so-called possibilistic integral.

The correspondents of these operators in the bivariate case can be defined on
$C([0,1]^{2})$ by the formulas%
\begin{multline*}
P_{n}(f)(x_{1},x_{2})=\sum_{k_{1}=0}^{n}\sum_{k_{2}=0}^{n}p_{n,k_{1}}%
(x_{1})p_{n,k_{2}}(x_{2})\\
\cdot\frac{\sup\{f(t_{1},t_{2})t_{1}^{k_{1}}(1-t_{1})^{n-k_{1}}t_{2}^{k_{2}%
}(1-t_{2})^{n-k_{2}}:t_{1},t_{2}\in\lbrack0,1]\}}{k_{1}^{k_{1}}n^{-n}%
(n-k_{1})^{n-k_{1}}k_{2}^{k_{2}}n^{-n}(n-k_{2})^{n-k_{2}}},
\end{multline*}
and%
\begin{multline*}
Q_{n}(f)(x_{1},x_{2})=\sum_{k_{1}=0}^{n}\sum_{k_{2}=0}^{n}p_{n,k_{1}}%
(x_{1})p_{n,k_{2}}(x_{2})\\
\cdot\sup\left\{  f(t_{1},t_{2}):t_{1}\in\left[  \frac{k_{1}}{n+1},\frac
{k_{1}+1}{n+1}\right]  ,t_{2}\in\left[  \frac{k_{2}}{n+1},\frac{k_{2}+1}%
{n+1}\right]  \right\} ,
\end{multline*}
respectively.

It is easy to show that both $P_{n}$ and $Q_{n}$ are monotone, unital and
sublinear operators. Notice that $Q_{n}$ is translatable, while $P_{n}$
is not.

Using the estimate included in the proof of Corollary 3.5 in \cite{Gal-3} (for
functions of one variable), one can easily show that
\begin{equation}
P_{n}(|\operatorname*{pr}\nolimits_{i}(\mathbf{t})-\operatorname*{pr}%
\nolimits_{i}(\mathbf{x})|)(\mathbf{x})\leq\frac{(1+\sqrt{2})\sqrt
{\operatorname*{pr}\nolimits_{i}(\mathbf{x})(1-\operatorname*{pr}%
\nolimits_{i}(\mathbf{x}))}+\sqrt{2}\sqrt{\operatorname*{pr}\nolimits_{i}%
(\mathbf{x})}}{\sqrt{n}}+\frac{1}{n}\,,\quad\label{estPn}%
\end{equation}
for $i\in\{1,2\},~n\in\mathbb{N},$ and all points $\mathbf{x}=(x_{1},x_{2})$
and $\mathbf{t}=(t_{1},t_{2})$ in $[0,1]^{2}.$

As a separating function on $[0,1]^{2}$ we choose the square distance,
\[
\gamma(\mathbf{x},\mathbf{y})=\left\Vert \mathbf{x}-\mathbf{y}\right\Vert
^{2}.
\]
We will show that%
\[
P_{n}(f)(x_{1},x_{2})\rightarrow f(x_{1},x_{2})\text{\quad uniformly on
}[0,1]^{2}%
\]
for all test functions \thinspace$1,~\pm\operatorname*{pr}\nolimits_{1}%
,~\pm\operatorname*{pr}\nolimits_{2},~\operatorname*{pr}\nolimits_{1}%
^{2}+\operatorname*{pr}\nolimits_{2}^{2}$ $($which will imply, via Theorem
$1$, that this convergence occurs for all nonnegative functions $f\in C([0,1]^{2}).$ As in
the case of Bernstein-Kantorovich-Choquet polynomial operators, only the
status of the test function $\operatorname*{pr}\nolimits_{1}^{2}%
+\operatorname*{pr}\nolimits_{2}^{2}$ needs attention. Or,%

\begin{multline*}
P_{n}(|\operatorname*{pr}\nolimits_{1}^{2}(\mathbf{t})+\operatorname*{pr}%
\nolimits_{2}^{2}(\mathbf{t})-\operatorname*{pr}\nolimits_{1}^{2}%
(\mathbf{x})+\operatorname*{pr}\nolimits_{2}^{2}(\mathbf{x})|)\\
\leq\sum_{i=1}^{2}P_{n}(|\operatorname*{pr}\nolimits_{i}(\mathbf{t}%
)+\operatorname*{pr}\nolimits_{i}(\mathbf{x})|\cdot|\operatorname*{pr}%
\nolimits_{i}(\mathbf{t})-\operatorname*{pr}\nolimits_{i}(\mathbf{x})|)\\
\leq2\sum_{i=1}^{2}P_{n}(|\operatorname*{pr}\nolimits_{i}(\mathbf{t}%
)-\operatorname*{pr}\nolimits_{i}(\mathbf{x})|)(x),
\end{multline*}
so that, according to (\ref{estPn}), we infer that $P_{n}(\operatorname*{pr}%
\nolimits_{1}^{2}+\operatorname*{pr}\nolimits_{2}^{2})\rightarrow
\operatorname*{pr}\nolimits_{1}^{2}+\operatorname*{pr}\nolimits_{2}^{2},$
uniformly on $[0,1]^{2}.$

The case of the operators $Q_{n}$ is similar. The fact that they verify the
hypotheses of Theorem 1 (for the same family of test functions) can be done as
above, by using instead the estimate included in the proof of Theorem 3.7 in
\cite{Gal-3} (for functions of one variable):
\[
Q_{n}(|\operatorname*{pr}\nolimits_{i}(\mathbf{t})-\operatorname*{pr}%
\nolimits_{i}(\mathbf{x})|)(\mathbf{x})\leq\frac{\sqrt{\operatorname*{pr}%
\nolimits_{i}(\mathbf{x})(1-\operatorname*{pr}\nolimits_{i}(\mathbf{x}))}%
}{\sqrt{n}}+\frac{2}{n+1}.
\]
Since $Q_{n}$ are translatable, the convergence of $Q_{n}(f)$ to $f$ holds for all $f\in C([0,1]^{2})$.

\medskip

\noindent\textbf{The max-product operators.} An important class of monotone,
unital and sublinear operators are the so-called max-product operators, whose
theory made the subject of the monograph \cite{BCG}. Denote $\bigvee_{j=0}%
^{m}=\max_{j=0,...,m}$ and $\Delta=\{(x_{1},x_{2});0\leq x_{1},x_{2}%
,x_{1}+x_{2}\leq1\}.$ The max-product Bernstein operators $T_{n}^{{}}%
:C(\Delta)\rightarrow C(\Delta)$ are defined by the formula
\[
T_{n}(f)(x_{1},x_{2})=\frac{\bigvee_{i=0}^{n}\bigvee_{j=0}^{n-i}{\binom{n}{i}%
}{\binom{n-i}{j}}x_{1}^{i}x_{2}^{j}(1-x_{1}-x_{2})^{n-i-j}f(i/n,j/n)}%
{\bigvee_{i=0}^{n}\bigvee_{j=0}^{n-i}{\binom{n}{i}}{\binom{n-i}{j}}x_{1}%
^{i}x_{2}^{j}(1-x_{1}-x_{2})^{n-i-j}}.
\]
As was shown in \cite{BCG}, pp. 139-140, these operators satisfy the estimate
\[
T_{n}(|\operatorname*{pr}\nolimits_{i}(\mathbf{t})-\operatorname*{pr}%
\nolimits_{i}(\mathbf{x})|)(\mathbf{x})\leq\frac{6}{\sqrt{n+1}},n\in
\mathbb{N},(x_{1},x_{2})\in\Delta,i=1,2.
\]
for $i\in\{1,2\},~n\in\mathbb{N},$ and all points $\mathbf{x}=(x_{1},x_{2})$
and $\mathbf{t}=(t_{1},t_{2})$ in $\Delta.$

The operators $T_{n}$ are sublinear unital and monotone (but not
translatable). Reasoning as in the previous example, one can infer from
Theorem \ref{thm1} that $T_{n}(f)(x_{1},x_{2})\rightarrow f(x_{1},x_{2})$
uniformly on $\Delta$ for all nonnegative functions $f\in C(\Delta)$.

As a consequence, for an arbitrary function $f\in C(\Delta)$ we have
\[
T_{n}(f+\left\Vert f\right\Vert _{\infty})-\left\Vert f\right\Vert _{\infty
}\rightarrow f\text{\quad uniformly on }\Delta.
\]

\noindent\textbf{The Gauss-Weierstrass-Choquet operators of two variables.}
The bivariate Gauss-Weierstrass-Choquet operators $W_{n,\mu}:C_{b}%
(\mathbb{R}^{2})\rightarrow C_{b}(\mathbb{R}^{2})$ are defined by the formula
\begin{multline*}
W_{n,\mu}(f)(x_{1},x_{2})\\
=\frac{(C)\int_{\mathbb{R}}(C)\int_{\mathbb{R}}f(s_{1},s_{2})e^{-n^{2}%
(x_{1}-s_{1})^{2}}\cdot e^{-n^{2}(x_{2}-s_{2})^{2}}d\mu(s_{1})d\mu(s_{2}%
)}{c(n,x_{1},\mu)c(n,x_{2},\mu)},
\end{multline*}
where $\mu=\sqrt{{\mathcal{L}}}$ and (according to the calculation in the
proof of Theorem 4.1 in \cite{Gal-6}) $c(n,x_{i},\mu)=(C)\int_{\mathbb{R}%
}e^{-n^{2}(x_{i}-s_{i})^{2}}d\mu(s_{i})=\sqrt{2/n}\cdot\Gamma(5/4)$, for
$i\in\{1,2\}$.

The fact that $W_{n,\mu}$ maps $C_{b}(\mathbb{R}^{2})$ into itself follows
from \cite{WK}, Theorem 11.13, p, 239.

Clearly, the operators $W_{n,\mu}$ are sublinear, monotone, unital but not
translatable.

Now, by using the estimate included in the proof of Theorem 4.1 in
\cite{Gal-6} (for functions of one variable), we infer that%
\[
W_{n,\mu}(|\operatorname*{pr}\nolimits_{i}(\mathbf{t})-\operatorname*{pr}%
\nolimits_{i}(\mathbf{x})|)(\mathbf{x})\leq\frac{4}{n}\quad\text{for }%
i\in\{1,2\},
\]
whenever $\mathbf{x}=(x_{1},x_{2})$ and $\mathbf{t}=(t_{1},t_{2})$ in
$\mathbb{R}^{2}.$ The Euclidean space $\mathbb{R}^{2}$ is locally compact
space and this also works for all equivalent metrics on it, in particular to
\[
d(\mathbf{x},\mathbf{y})=|x_{1}-y_{1}|+|x_{2}-y_{2}|,
\]
for all $\mathbf{x}=(x_{1},x_{2})$ and $\mathbf{y}=(y_{1},y_{2})$ in
$\mathbb{R}^{2}$. Taking into account Theorem \ref{thm4} (for $p=1$), one can
easily show that $W_{n,\mu}(f)\rightarrow f,$ uniformly on the compact subsets
of $\mathbb{R}^{2}$, for every nonnegative $f\in C_{b}(\mathbb{R}^{2}).$

\section{Appendix: Generalities on Choquet's integral}

Very interesting and the integral associated to it. Full details are to be
found in the books of D. Denneberg \cite{Denn}, M. Grabisch \cite{Gr2016} and
Z. Wang and G. J. Klir \cite{WK}.

Let $(X,\mathcal{A})$\ be an arbitrarily fixed measurable space, consisting of
a nonempty abstract set $X$ and a $\sigma$-algebra ${\mathcal{A}}$ of subsets
of $X.$

\begin{definition}
\label{def1} A set function $\mu:{\mathcal{A}}\rightarrow\lbrack0,1]$ is
called a capacity if it verifies the following two conditions:

$(a)$ $\mu(\emptyset)=0$ and $\mu(X)=1;$

$(b)~\mu(A)\leq\mu(B)$ for all $A,B\in{\mathcal{A}}$, with $A\subset B$
\emph{(}monotonicity\emph{)}.
\end{definition}

An important class of capacities is that of probability measures (that is, the
capacities playing the property of $\sigma$-additivity). Probability
distortions represents a major source of nonadditive capacities. Technically,
one start with a probability measure $P:\mathcal{A\rightarrow}[0,1]$ and
applies to it a distortion $u:[0,1]\rightarrow\lbrack0,1],$ that is, a
nondecreasing and continuous function such that $u(0)=0$ and $u(1)=1;$for
example, one may chose $u(t)=t^{a}$ with $\alpha>0.$The \emph{distorted
probability} $\mu=u(P)$ is a capacity with the remarkable property of being
continuous by descending sequences, that is,%
\[
\lim_{n\rightarrow\infty}\mu(A_{n})=\mu\left(
{\displaystyle\bigcap_{n=1}^{\infty}}
A_{n}\right)
\]
for every nonincreasing sequence $(A_{n})_{n}$ of sets in $\mathcal{A}.$ Upper
continuity of a capacity is a generalization of countable additivity of an
additive measure. Indeed, if $\mu$ is an additive capacity, then upper
continuity is the same with countable additivity. When the distortion $u$ is
concave (for example, when $u(t)=t^{a}$ with $0<\alpha<1),$ then $\mu$ is also
\emph{submodular} in the sense that%
\[
\mu(A\cup B)+\mu(A\cap B)\leq\mu(A)+\mu(B)\text{\quad for all }A,B\in
\mathcal{A}.
\]

The next concept of integrability with respect to a capacity refers to the
whole class of random variables, that is, to all functions $f:X\rightarrow
\mathbb{R}$ such that $f^{-1}(A)\in{\mathcal{A}}$ for every Borel subset $A$
of $\mathbb{R}$.

\begin{definition}
\label{def2}The Choquet integral of a random variable $f$ with respect to the
capacity $\mu$ is defined as the sum of two Riemann improper integrals,
\[
(C)\int_{X}fd\mu=\int_{0}^{+\infty}\mu\left(  \{x\in X:f(x)\geq t\}\right)
dt+\int_{-\infty}^{0}\left[  \mu\left(  \{x\in X:f(x)\geq t\}\right)
-1\right]  dt,
\]
Accordingly, $f$ is said to be Choquet integrable if both integrals above are finite.
\end{definition}

If $f\geq0$, then the last integral in the formula appearing in Definition
\ref{def2} is 0.

The inequality sign $\geq$ in the above two integrands can be replaced by $>;$
see \cite{WK}, Theorem 11.1,\emph{ }p. 226.

Every bounded random variable is Choquet integrable. The Choquet integral
coincides with the Lebesgue integral when the underlying set function $\mu$ is
a $\sigma$-additive measure.

As usually, a function $f$ is said to be Choquet integrable on a set
$A\in\mathcal{A}$ if $f\chi_{A}$ is integrable in the sense of Definition
\ref{def2}. We denote%
\[
(C)\int_{A}fd\mu=(C)\int_{X}f\chi_{A}d\mu.
\]
We next summarize some basic properties of the Choquet integral.

\begin{remark}
\label{remApp}$(a)$ If $\mu:{\mathcal{A}}\rightarrow\lbrack0,1]$ is a
capacity, then the associated Choquet integral is a functional on the space of
all bounded random variables such that:%
\begin{gather*}
f\geq0\text{ implies }(C)\int_{A}fd\mu\geq0\text{ \quad\emph{(}%
positivity\emph{)}}\\
f\leq g\text{ implies }\left(  C\right)  \int_{A}fd\mu\leq\left(  C\right)
\int_{A}gd\mu\text{ \quad\emph{(}monotonicity\emph{)}}\\
\left(  C\right)  \int_{A}afd\mu=a\cdot\left(  \left(  C\right)  \int_{A}%
fd\mu\right)  \text{ for }a\geq0\text{ \quad\emph{(}positive\emph{
}homogeneity\emph{)}}\\
\left(  C\right)  \int_{A}1\cdot d\mu(t)=\mu(A)\text{\quad\emph{(}%
calibration\emph{)}};
\end{gather*}
see \emph{\cite{Denn},} Proposition \emph{5.1} $(ii)$, p. \emph{64}, for a
proof of the property of positive\emph{ }homogeneity.

$(b)$ In general, the Choquet integral is not additive but, if the bounded
random variables $f$ and $g$ are comonotonic,\emph{ }then%
\[
\left(  C\right)  \int_{A}(f+g)d\mu=\left(  C\right)  \int_{A}fd\mu+\left(
\operatorname*{Ch}\right)  \int_{A}gd\mu.
\]
This is usually referred to as the property of comonotonic additivity and was
first noticed by Delacherie \cite{Del1970}. An immediate consequence is the
property of translation invariance,
\[
\left(  C\right)  \int_{A}(f+c)d\mu=\left(  C\right)  \int_{A}fd\mu+c\cdot
\mu(A)
\]
for all $c\in\mathbb{R}$ and all \hspace{0in}bounded random variables $f.$ For
details, see \emph{\cite{Denn}}, Proposition \emph{5.1, }$(vi)$, p. \emph{65.}

$(c)$ If $\mu$ is a lower continuous capacity, then the Choquet integral is
lower continuous in the sense that%
\[
\lim_{n\rightarrow\infty}\left(  \left(  C\right)  \int_{A}f_{n}d\mu\right)
=\left(  C\right)  \int_{A}fd\mu
\]
whenever $(f_{n})_{n}$ is a nondecreasing sequence of bounded random variables
that converges pointwise to the bounded variable $f.$ See \emph{\cite{Denn}},
Theorem $8.1$, p. $94.$

$(d)$ Suppose that $\mu$ is a submodular capacity. Then\ the associated
Choquet integral is a subadditive functional, that is,
\[
\left(  C\right)  \int_{A}(f+g)d\mu\leq\left(  C\right)  \int_{A}fd\mu+\left(
C\right)  \int_{A}gd\mu
\]
for all bounded random variables $f$ and $g.$ See \emph{\cite{Denn}},
Corollary \emph{6.4}, p. \emph{78.} and Corollary\emph{ 13.4, }p.\emph{ 161.
}It is also a submodular functional in the sense that
\[
\left(  C\right)  \int_{A}\sup\left\{  f,g\right\}  d\mu+\left(  C\right)
\int_{A}\inf\{f,g\}d\mu\leq\left(  C\right)  \int_{A}fd\mu+(C)\int_{A}gd\mu
\]
for all bounded random variables $f$ and $g.$ See \cite{CMMM2012}, Theorem
$13$, $(c)$.
\end{remark}

\end{document}